# Parameter Estimation of Differential Equation Model based on Optimal Weight Choice Method

Jun Wang[a], Xianglei Li [a] , Xianghu Li[a]

[a]*School of Mathematics and Statistics, Changchun University of Technology, Changchun 130012, China*

(Changchun University of Technology,School of Mathematics and Statistics,Changchun Jilin 130012)

**Abstract**
Differential equations are important tools to portray dynamic problems, and are widely used in finance, engineering and biology. Here, multiple dynamic differential models were built innovatively, and discretized with the Runge-Kutta method. The the model parameters were estimated. The models were averaged using the Optimal weight selection method, and the consistency of such parameter estimation was verified. Numerical simulations were also conducted, and the simulated results outperformed ordinary linear models. Finally, the differential averaging model built here was used to empirically analyze the Shanghai Index 300. This method integrated the fluctuation features of multiple Shanghai Composite Index fitting models, and yielded good analytical results. This study provides a methodological reference for analysis of stock market situations, and offers research clues for the parameter estimation of differential equations.

**Key words** differential equation · parameter estimation · Optimal Weight Choice

## 1  Introduction

Ordinary differential equations are widely used to describe the complex dynamic systems of finance, management science, ecology and other scientific fields. Especially in finance statistics, ordinary differential equations are commonly used to analyze the trend of stock indices, and the increasing rate variation of national GDP with time. There are many excellent models and methods for parameter estimation and variable selection of static systems. With consideration into the merits of these models, the model averaging method can be used. Currently, model averaging has developed mainly to two branches: frequency model averaging, and Bayesian model averaging. The prior distribution of the model shall be considered during Bayesian model averaging, which may be difficult in some situations. Thus, frequency model averaging has been greeted by more experts and researchers owing to its practicability and convenience. Frequency model averaging, also called the model composite method, assigns each submodel of a composite model with a specific weight, and is prospective in finance, medicine and other fields. This study is aimed to efficiently expand the frequency model averaging to parameter estimation and variable selection of ordinary differential equations. This model will largely improve the accuracy in parameter estimation of ordinary differential equations.

Bates and Granger were the first to find that the statistical deduction of traditional model selection methods is based on a given linear model, but the specific information reflected by other models may be ignored. Hjort and Claeskens found traditional model selection methods did not consider the uncertainty due to the way of selection, so the final model-based deduction or prediction may contain much error [1]. Nevertheless, the appearance of model averaging solved all the defects of traditional model selection methods, and can better explain the existing data [2]. This method can be traced back to 1960s, and its core idea is to assign multiple models with appropriate weights according to known data. It provides an equilibrium mechanism: the uncertainty of each candidate model is balanced by integrating the parameter estimations of multiple models. Roberts combined the distributive characteristics of other models and weight-averaged the posteriori distributions of two types of models [3]. This method was similar to Bayesian model averaging, and the uncertainty of model selection was also considered. Buckland and Burnham proposed Smoothed Akaike information criterion (S-AIC) and smoothed Bayesian information criterion (S-BIC), and thereby established the most basic model averaging method [4], which is commonly used owing to its easy operations. In the 21th century, Yuan and Yang put forward the concrete steps of model selection, and achieved adaptive regression in a mixed way [5]. Hansen proposed Mallows model averaging (MMA), in which the weights of candidate models were determined from Mallows criterion of minimization [6]. This method has largely contributed to the development of model averaging. The Mallows criterion is nearly equivalent to the square errors of models. Thus, the Mallows criterion of minimization is equivalent to minimizing the mean square error of estimation. On this basis, Liang H proposed the Optimal model averaging for weight selection [7]. Compared with MMA, this method is characteristic of small samples, and is nearly the optimal relative to other model selection methods. Later, Hansen and Racine designed Jackknife model averaging (JMA) [8], which determines weights via the Jackknife criterion of minimization like other estimation methods. They also validated that JMA is asymptotically optimal with the presence of heteroskedasticity in models. As for the nesting models of model averaging, Zhang and Liu studied the distribution of parameters estimated by the least squares method under fixed parameters and model nesting, and clarified the asymptotic distribution of estimations from MMA and JMA [9]. The confidence intervals of parameter estimations were also provided. Moreover, when the candidate models involved a real model, it was proved that under the Mallows criterion, the model with deficient goodness-of-fit will have a weight converging to 0 according to the probability. The weights of the over-fitting model and the real model were asymptotically random, and the estimated values of parameters converged to a non-standard distribution. Based on the findings from Hansen et al., Xinyu Zhang et al. developed prediction model averaging (PMA), which was mainly applied to high-dimension data models and was proved to estimate model parameters that are asymptotically the optimal [10]. Xinmin Li et al. further proposed a weight selection method based on generalized cross-validation approximation [11]. While the least possible squared error is achieved, the estimations from least square model averaging are proved to be asymptotically optimal. Especially, the optimum is achieved under a disperse weight set and a continuous weight set. Compared with the existing methods based on the Mallows criterion, the asymptotic optimum conditions required by this method are more reasonable.

In all, the model averaging methods proposed recently all validate the asymptotic optimum of model averaging estimation. However, the existing model averaging methods are almost all based on general linear models. Herein, together with ordinary differential equations, the excellent quality of practical problems was portrayed, and model averaging methods were used to estimate model parameters. As for model averaging estimation in the case of small samples, Optimal weight selection method performs well in practice. Thus, Optimal model averaging was used for model averaging estimation of differential equations.

This study is organized as follows: Section 2 presents modeling and the parameter estimation of differential equations. Section 3 is about numerical simulation. Section 4 shows empirical analysis with Shanghai Composite Index data. Section 5 lists the conclusions. Consistence certification is provided in the supporting materials.

## 2. Modeling and parameter estimation

### 2.1 Fourth-order Runge-Kutta method for differential model

In a general differential equation

$$\frac{dY}{dt} = F(t, Y, \alpha),$$

where $Y(t) = (Y_1(t), Y_2(t), \cdots Y_p(t))^T$ is a p-dimension vector, $\alpha = (\alpha_1, \alpha_2, \cdots \alpha_q)^T$ is parameter vector, and $F(\cdot)$ is a p-dimension function. Currently, many methods can be used to solve dynamic differential equations, such as $Euler$ method and $Heun$ method. Though their physical meaning in solving differential equations is clear and they are easy-to-program, their computational accuracy is not high theoretically. If the step length $h$ is too small, the number of interval segments certainly will increase. Owing to limitations by computer capability and valid digit number, error accumulation is unavoidable, so the results decrease computational accuracy instead. To improve precision, we adopted four-order Runge-Kutta (short for RK4) to discretize differential equations. RK4 can be used to simulate and compute according to known differential equations and initial values while ignoring the complexity of solving, and thus becomes a commonly-used method in computing differential equations.

The Runge-Kutta method is derived from $Taylor$ expansion. The commonly-used fourth-order Runge-Kutta method (RK4) is elaborated below.

In a first-order ordinary differential equation

$$\frac{dy}{dt} = f(t, y),$$

$$\begin{cases} y_{i+1} = y_i + \dfrac{h}{6}(k_1 + 2k_2 + 2k_3 + k_4) \\ k_1 = f(t_i, y_i) \\ k_2 = f(t_i + h/2, y_i + \dfrac{h}{2}k_1) \\ k_3 = f(t_i + h/2, y_i + \dfrac{h}{2}k_2) \\ k_4 = f(t_i + h, y_i + hk_3) \end{cases} \quad (1)$$

$h$ is the step length and its value can be set freely; its truncation error is $O(h^5)$. $y_{i+1}$ is computed

from $y_i$, and the function value is computed four times with each step forward. In other words, based on a given function $f(t_i, y_i)$, $y_1, y_2, ..., y_n$ can be determined from the initial value $y_0$. For the RK4 of a high-order differential equation, each component is subjected to RK4.

Generally, F can be expressed in any linear or nonlinear function. However, the nonlinearity of F calls for priori information and high computational cost. Thus, nonlinear ODE models are only applicable to small-scale problems that contain several to tens of explaining variables. Many differential models are based on linear ODE owing to its simpleness and practicability. However, we realize the dynamic system of any problem may exist in complex pattern, which may not be fully covered by a linear model. For simplification of narration, we write the ODE model in a simple form:

$$\frac{dY}{dt} = AY(t) \tag{2}$$

where $Y(t) = (Y_1(t), Y_2(t), \cdots Y_d(t))^T$ is an n-dimension explanatory variable, with a coefficient matrix $A = (\alpha_{ij}) \in \Re^{p \times d}$. In the subsequent texts, we will simplify time variables, such as $Y(t) = Y_t$. The dynamic differential model (2) is discretized using RK4 in Eq. (1) as follows:

$$Y(t+1) - Y(t) \approx \frac{h}{6} A\left(Y_t + 4Y_{t+h/2} + Y_{t+h}\right) \tag{3}$$

As there is random error $\varepsilon_i$ between the observation data $y$ and $Y$:

$$y = Y + \varepsilon, \qquad \varepsilon \sim N(0, \sigma^2 I_\alpha).$$

from Eq. (3), we have a regression model

$$\Delta y_t = \frac{h}{6} A\left(Y_t + 4Y_{t+h/2} + Y_{t+h}\right) + h\varepsilon, \tag{4}$$

where $E(\varepsilon | x) = E(\varepsilon_t + 4\varepsilon_{t+h/2} + \varepsilon_{t+h}) = 0$. Let function $D_k$ be the sum of squared errors from the k-th row of components in Eq. (4),

$$D_k = \sum_{i=1}^{n}\left(\Delta y_{k,t_i} - \frac{h}{6}\alpha_{k\bullet}\left(Y_{t_i} + 4Y_{t_i+h/2} + Y_{t_i+h}\right)\right),$$

The partial derivative of $D_k$ relative to $\alpha_{k\bullet}$ is solved: $\Delta D_k = 0$. In this way, the model parameter $\alpha_{k\bullet}$ can be determined

$$\hat{\alpha}_{k\bullet}^T = \frac{6\sum_{i=1}^{n}\left((\Delta y_{k,t_i} - \Delta \bar{y}_k)(Y_{t_i} + 4Y_{t_i+h/2} + Y_{t_i+h} - 6\bar{Y})\right)}{h\sum_{i=1}^{n}\left((Y_{t_i} + 4Y_{t_i+h/2} + Y_{t_i+h} - 6\bar{Y})^T(Y_{t_i} + 4Y_{t_i+h/2} + Y_{t_i+h} - 6\bar{Y})\right)}.$$

## 2.2 Optimal weight selection based on frequency model averaging

Without missing generality, we only analyze k rows (k=1,2,......, p) in Eq. (4). For simplification, we leave out the mark k. Frequency model averaging is used to discretize differential equation:

$$\Delta y_t = (x_t + 4x_{t+h/2} + x_{t+h})\frac{h}{6}\beta + (z_t + 4z_{t+h/2} + z_{t+h})\frac{h}{6}\gamma + \frac{h}{6}(\varepsilon_t + 4\varepsilon_{t+h/2} + \varepsilon_{t+h}), \quad (5)$$

where $\Delta y_t$ is the first-order differential variable of the dependent variable $y$. After transposition of $Y$, n groups of data are collected to constitute a regression variable data matrix $(X,Z)$, $X = x_t + 4x_{t+h/2} + x_{t+h}$ is an $n \times k$ matrix, $Z = z_t + 4z_{t+h/2} + z_{t+h}$ is an $n \times m$ matrix, $d = k + m$. Moreover, $x_t$ is the main variable among the regression variables (meaning all variables that must be included in the model); $z_t$ is the secondary variable (meaning the variables that are nonessential in the model). Such setting will not miss universality, as the matrix $X$ may not contain any regression factor except for an intercept term and may even be empty. $\beta$ and $\gamma$ are the parameters to be estimated: $\alpha_{k\bullet} = (\beta, \gamma)^T$. As $Z$ contains $m$ optimal auxiliary regression variables, there are at most $N = 2^m$ expanded models available.

If $\gamma = 0$, let

$$\hat{\beta}_r = \frac{6}{h}(X'X)^{-1}X'\Delta y,$$

If $\gamma \neq 0$, let

$$\hat{\beta}_u = \hat{\beta}_r - \frac{6}{h}(X'X)^{-1}X'Z(Z'MZ)^{-1}Z'M\Delta y,$$

Corresponding with $\hat{\gamma} = \frac{6}{h}(Z'MZ)^{-1}Z'M\Delta y$, where $M = I_n - X(X'X)^{-1}X'$.

Now, let $\theta = \frac{h}{6}(Z'MZ)^{1/2}\gamma$, and $\hat{\theta} = \frac{h}{6}(Z'MZ)^{1/2}\hat{\gamma}$. We note $\hat{\theta} \sim N(\theta, \sigma^2 I_m)$, where

$$\sigma^2 = \|\Delta y - X\frac{h}{6}\hat{\beta}_u - Z\frac{h}{6}\hat{\gamma}\|^2 / n.$$

Then let $Q = (X'X)^{-1}X'Z(Z'MZ)^{-1/2}$, for easy expression that $\frac{h}{6}\hat{\beta}_u = \frac{h}{6}\hat{\beta}_r - Q\hat{\theta}$. The unestimated parameter $\beta$ of the ith model can be expressed as $\hat{\beta}_{(i)} = \hat{\beta}_r - \frac{6}{h}QW_i\hat{\theta}$, where

$W_i = I_m - P_i$, $P_i = (Z'MZ)^{-1/2}S_i\{S'_i(Z'MZ)^{-1}S_i\}^{-1}S'_i(Z'MZ)^{-1/2}$ is the selection matrix.

Based on the above conditions, the model averaging parameter estimation expression can be shown below：

$$\hat{\beta}_f = \sum_{i=1}^{N} \lambda_i \hat{\beta}_i, \quad \hat{\gamma}_f = \sum_{i=1}^{N} \lambda_i \hat{\gamma}_i, \tag{6}$$

where $\lambda_i \geq 0$ is the weight to be assigned, and $\sum_{i=1}^{N} \lambda_i = 1$. Let $\sigma_i^2 = \| \Delta y - \frac{h}{6} X \hat{\beta}_{(i)} - \frac{h}{6} Z \hat{\gamma}_i \|^2 / n$ be the mean square error of the $i$-th model. Notably, $\sigma_i^2$ can also be expressed as $\sigma_i^2 = (n - k - m)\sigma^2 / n + \hat{\theta}' P_i \hat{\theta} / n$. When this equation is substituted to the AIC expression of the $i$-th model, we find the value of AIC is decided only by $\hat{\theta}$ and $\sigma^2$. Hence, it can be stated that the weight coefficient $\lambda_i = \lambda_i(\hat{\theta}, \sigma^2)$ is also only related to $\hat{\theta}$ and $\sigma^2$. Next, $\lambda_i$ is estimated. The weight selection method is set according to the target value of MSE minimization, especially in the case of a finite sample size. The MSE of $\hat{\beta}_f$ is $E\{(\hat{\beta}_f - \beta)(\hat{\beta}_f - \beta)'\}$, diagonal elements of which are the MSEs of estimation of $\beta$. The trace of the MSE matrix is equal to the expected squared error loss function $E \| \hat{\beta}_f - \beta \|^2$, and is also called the estimator risk under squared error loss, or weak MSE precision tolerance [12]. The MSE of $\hat{\beta}_f$ can be expressed as follows [7]:

$$MSE(\hat{\beta}_f) = \hat{\sigma}^2 (X'X)^{-1} - \hat{\sigma}^2 QQ' + \{Q(I_m - W)\hat{\theta}\}^{\otimes 2} + \psi(\hat{\theta}, \hat{\sigma}^2) + \{\psi(\hat{\theta}, \hat{\sigma}^2)\}'.$$

The trace of $MSE(\hat{\beta}_f)$ is:

$$\hat{R}(\hat{\beta}_f) = \hat{\sigma}^2 tr(X'X)^{-1} - \hat{\sigma}^2 tr(QQ')$$
$$+ \{\hat{\theta}'(I_m - W)Q'\}^{\otimes 2} + 2\hat{\sigma}^2 tr\{\psi_1(\hat{\theta}, \hat{\sigma}^2)\}, \tag{7}$$

where $\psi_1(\hat{\theta}, \hat{\sigma}^2) = QW + Q \sum_{i=1}^{N} \{\partial \lambda_i(\hat{\theta}, \hat{\sigma}^2) / \partial \hat{\theta}\} \hat{\theta}' W_i$.

For frequency model averaging, the unified weight expression is

$$\lambda_i(\hat{\theta}, \hat{\sigma}^2) = \frac{a^{q_i}(n - q_i)^b (\hat{\sigma}_i^2)^c}{\sum_{j=1}^{N} a^{q_j}(n - q_j)^b (\hat{\sigma}_j^2)^c}, \tag{8}$$

where $q_i$ is the number of regression variables in each model, and $a > 0, b > 0, c < 0$ are constants. Let

$$\partial \lambda_i(\hat{\theta}, \hat{\sigma}^2) / \partial \hat{\theta} = \frac{2}{n} c \lambda_i(\hat{\theta}, \hat{\sigma}^2) \{(\hat{\sigma}_i^2)^{-1}(I_m - W_i)$$
$$- \sum_{j=1}^{N} \lambda_j(\hat{\theta}, \hat{\sigma}^2)(\hat{\sigma}_j^2)^{-1}(I_m - W_j)\} \hat{\theta}. \tag{9}$$

Eqs. (7), (8) and (9) are combined to minimize $\hat{R}(\hat{\beta}_f)$, so the weights of coefficients in each model can be determined. On this basis, the values of coefficients $\hat{\beta}_f$ and $\hat{\gamma}_f$ can be estimated.

### 2.3. Consistence of Optimal estimation based on differential model

The discussion in section 2.2 is focused on Optimal estimation with a finite sample size. Here, the consistency of estimation with large sample size will be analyzed. The linear model converted from the Runge-Kutta process is

$$\Delta \hat{y}_t = \hat{\mu}_t = \frac{h}{6}\hat{\beta}\left(x_t + 4x_{t+h/2} + x_{t+h}\right) + \frac{h}{6}\left(\varepsilon_t + \varepsilon_{t+h/2} + \varepsilon_{t+h}\right),$$

The square error loss function of estimation is $Ln(\lambda(a,b,c)) = \|\hat{\mu}_f(\lambda(a,b,c)) - \mu\|^2$, and the risk loss function under squared errors is $R_n(\lambda(a,b,c)) = E\{Ln(\lambda(a,b,c))\}$. Let the space be $D = \{(a,b,c) \mid a > 0,\ b \geq 0, -\bar{c} \leq c \leq 0\}$, where $\bar{c}$ is a constant larger than 0. Let the set $U$ refer to an unbiased model, or namely, a real model is nested as a special model. Besides, the following subset of $D$ is considered

$$D_0 = \left\{(a,b,c) \in D \mid \sum_{\tau \in U} \lambda_\tau(a,b,c) \leq 1 - \rho \right\},$$

where $\rho$ is a constant in the range $(0,1]$. In $D_0$, the weights at $\rho = 0$ assigned to the biased model have a nonzero lower limit during the model averaging. In this row, all cases are unbiased except for the submodel during the model averaging. The loss of generality due to such restriction is very small. It can be imagined that the model averaging method is more likely to contain some biased models rather than unbiased models. Based on this framework, we define a non-random weight set

$$W = \left\{w = (w_1, \ldots, w_N)' \mid w_i \geq 0,\ \sum_{i=1}^{N} w_i = 1,\ \sum_{\tau \in U} w_\tau \leq 1 - \rho \right\}.$$

Let $\zeta_n = \max_{1 \leq i \leq N} E\|\hat{\mu}_i - \mu\|^2$ be the maximum risk amount in a single sub-model, and $\zeta_n = \inf_{w \in W} R_n(w)$. Let $(\hat{a}, \hat{b}, \hat{c})$ be the value of $(a,b,c)$ after $\hat{R}_a(\hat{\mu}_f(\lambda(a,b,c)))$ is minimized, and $(\hat{a}, \hat{b}, \hat{c}) \in D_0$. Then we have the following theorem.

**Theorem 1:** When $n \to \infty$, and conditions $\mu'\mu = O(n)$ and $\zeta_n^{-2} \zeta_n \to 0$ are met, we have

$$\frac{Ln(\lambda(\hat{a},\hat{b},\hat{c}))}{\inf_{(a,b,c) \in D_0} Ln(\lambda(a,b,c))} \xrightarrow{p} 1.$$

This theorem indicates the consistency of Optimal parameter estimation is true in the case of large sample sizes. The proving is similar to the case of static linear models and is shown in the appendix.

## 3. Simulation

In this section, data are generated from the following model

$$\frac{dY_i}{dt} = \alpha_{i\bullet} Y(t) + \varepsilon,$$

where $Y_i$ are 4 inter-independent series of Y in specific distribution, and $\varepsilon \sim N(0,1)$. Let $Y_1$ and $Y_2$ be the main variables in the model, and $Y_3$ and $Y_4$ be auxiliary variables. In this way, the weights of totally $N = 2^2$ (or 4) models shall be computed. The sample size is between $n = 50$ and 300. The error term $\varepsilon$ is independent relative to all $Y_i$. In the simulation, we focus on how the first two main variables change with the convergence of sample size variation during Optimal model averaging weight selection, and on comparing the overall mean square errors between differential equations and linear models.

### 3.1 First group of simulation

Under the above hypotheses, let the real value of the estimation coefficient be $\alpha = (\alpha_1, \alpha_2, \alpha_3, \alpha_4) = (0.1, 0.2, 0.3, 0.4)$, and the generated data obey $Y_1 \sim N(1, 0.2)$, $Y_2 \sim N(2, 0.3)$, $Y_3 \sim N(3, 0.4)$ and $Y_4 \sim N(4, 0.3)$. The the simulated results with sample sizes $n = 200$ and $n = 300$ will be shown.

Table 1  Analysis of simulation with $n = 200$

| Real value | Simulated value | Deviation | Mean square error |
|---|---|---|---|
| 0.1 | 0.0934 | -0.0066 | 0.1425 |
| 0.2 | 0.2050 | 0.0050 | 0.0712 |

Table 2  Analysis of simulation with $n = 300$

| Real value | Simulated value | Deviation | Mean square error |
|---|---|---|---|
| 0.1 | 0.0989 | 0.0011 | 0.1265 |
| 0.2 | 0.1932 | -0.0068 | 0.0487 |

When at $n = 200$ and $n = 300$, the simulated value is very close to the real value, and the deviation and mean square error (MSE) of computation are both small (Tables 1 and 2). Moreover, at the sample size $n = 300$, the deviation and MSE are both smaller than at sample size $n = 200$.

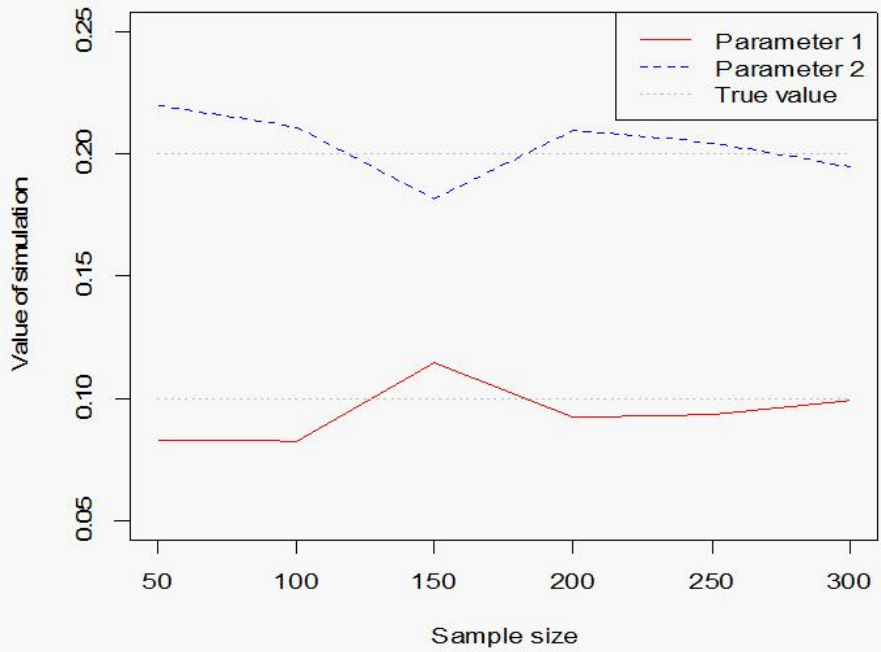

**Fig. 1** The trend diagram of simulation

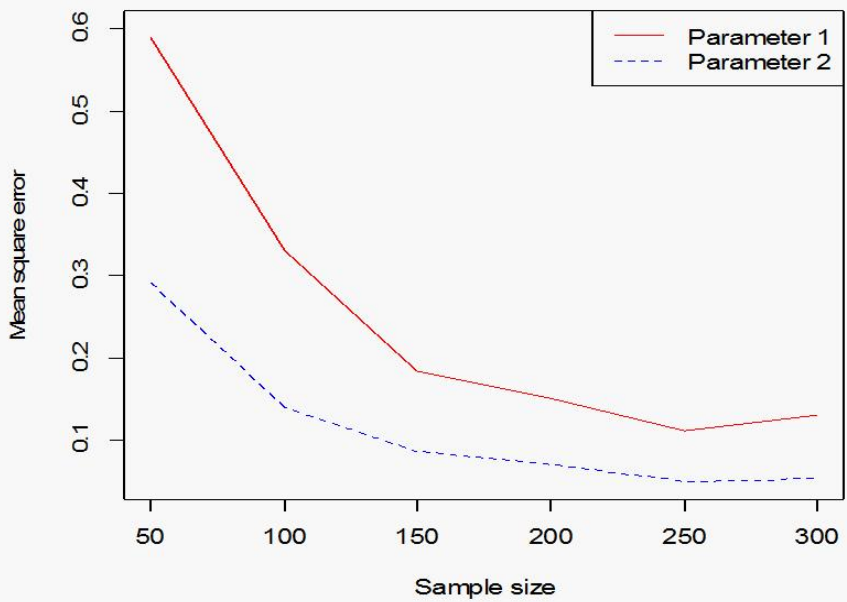

**Fig. 2** The trend of the mean square error of the parameter

The curves of simulation and the MSEs of parameter estimation show that as the sample size rises, the MSEs of parameter estimation gradually decrease, and the values of parameter estimation basically converge to the real value. Thus, the consistency of parameter estimation is validated.

Then the overall MSEs of differential equations and linear models are compared during simulation. In the small sample size 200, let $MSE1$ be the MSEs of differential equations, and $MSE2$ be the MSEs of linear models. The comparison is listed in Table 3.

Table 3 Comparison of simulated mean square error

| Sample size of $n$ | 40 | 80 | 120 | 160 | 200 |
|---|---|---|---|---|---|
| $MSE1$ | 0.2099 | 0.6738 | 0.6032 | 0.7047 | 0.5928 |
| $MSE2$ | 0.8271 | 1.2088 | 0.8617 | 0.9802 | 1.0039 |

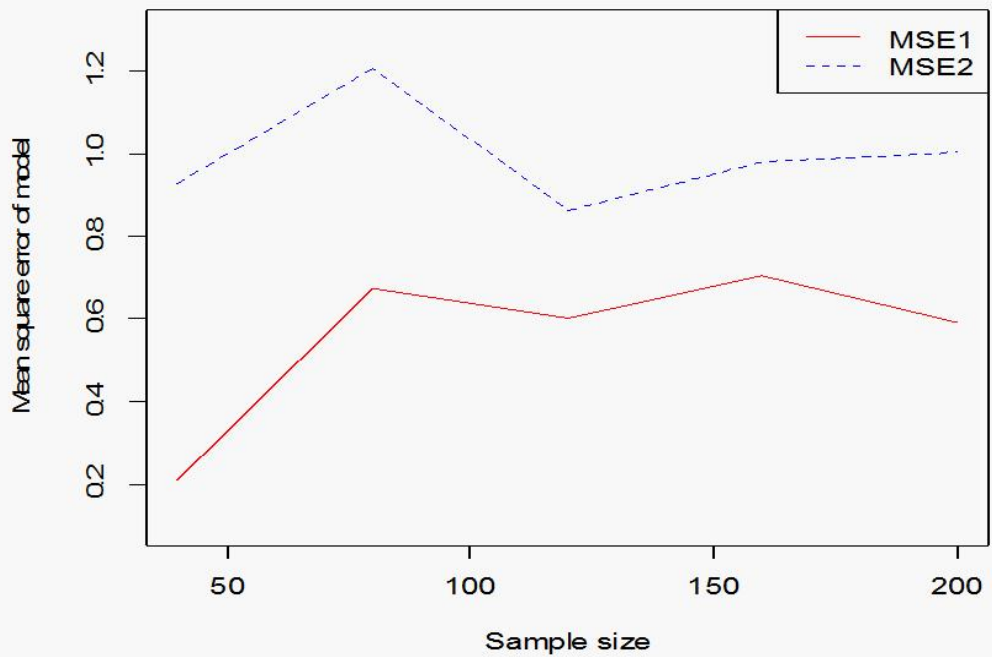

Fig. 3  Comparison of simulated mean square error

At the same sample size, the MSEs of differential equations are smaller than those of linear models (Table 3). The amplitude of fluctuation of MSEs from the differential equations are smaller compared with the linear model (Fig. 3). Thus, during modeling, the differential equation averaging model theoretically has smaller MSEs than the linear average model, and thus fits better.

**3.2 Second group of simulation**

Let the real value of the estimation coefficient be $\alpha = (\alpha_1, \alpha_2, \alpha_3, \alpha_4) = (0.4, 0.5, 0.6, 0.7)$, and the generated data be $Y_1 \sim N(1.5, 0.3)$, $Y_2 \sim N(2.5, 0.4)$, $Y_3 \sim N(3.5, 0.5)$ and $Y_4 \sim N(4.5, 0.6)$. The the simulated results with sample sizes $n = 200$ and $n = 300$ will be shown.

Table 4　Analysis of simulation with $n=200$

| Real value | Simulated value | Deviation | Mean square error |
|---|---|---|---|
| 0.4 | 0.4143 | 0.0143 | 0.1081 |
| 0.5 | 0.4975 | -0.0025 | 0.0597 |

Table 5　Analysis of simulation with $n=300$

| Real value | Simulated value | Deviation | Mean square error |
|---|---|---|---|
| 0.4 | 0.4030 | 0.0030 | 0.0868 |
| 0.5 | 5153 | 0.0153 | 0.0408 |

When at $n=200$ and $n=300$, the simulated value is very close to the real value, and the deviation and MSE of computation are both small (Tables 4 and 5). Moreover, at the sample size $n=300$, the MSEs are both smaller than at sample size $n=200$.

Then the overall MSEs of differential equations and linear models are compared during simulation. In the small sample size 200, let $MSE1$ be the MSEs of differential equations, and $MSE2$ be the MSEs of linear models. The comparison is listed in Table 6.

Table 6　Comparison of simulation mean square errors

| Sample size of $n$ | 40 | 80 | 120 | 160 | 200 |
|---|---|---|---|---|---|
| $MSE1$ | 0.3082 | 0.4231 | 0.2965 | 0.3586 | 0.3132 |
| $MSE2$ | 0.7521 | 0.6746 | 0.5256 | 0.6893 | 0.4937 |

At the same sample size, the MSEs of differential equations are smaller than those of linear models (Table 6). The fluctuation amplitude of MSEs from the differential equations are smaller compared with the linear models (Fig. 3). Thus, during modeling, the differential model averaging theoretically has smaller MSEs than the linear model averaging, and thus fits better.

Based on the above two groups of simulations, as for the same group of data, the differential equation model outperforms traditional linear average models in terms of fitting.

## 4. Empirical analysis

Recently, the Shanghai Composite Index is used in the Chinese stock market to reflect changes in the stock market, and is a key data source used by experts and researchers in various fields. Along with economic globalization, the Chinese stock market has developed well in recent years. Shanghai Stock Exchange has been founded for 30 years. The trends of Shanghai Composite Index reflect the variations of stock market prices to large extent. Hence, the fluctuation and trends of this index are studied and predicted by finance and economy researchers using various models. Herein, a stock index dynamic differential model was built first, and then its parameters were estimated using the Optimal weight selection method. Finally, the advantage of this differential model in portraying Shanghai Composite Index was validated by comparing with traditional models.

**4.1 Data source**

Totally 57 groups of data were collected from the daily data of Shanghai Composite Index between January 2016 and October 2020 on Straight Flush. The dependent variable $y_0$ is the closing price of Shanghai Composite Index on the last trading day of each month, and the independent variables $y_1 \sim y_8$ are the amplitude, total lots, difference exponential average (DEA), price movement, DIFF, rising amplitude, and 5-day average. The definitions of the variables are listed in Table 7.

Table 7  Introduction of variables

| Variable | Description |
| --- | --- |
| Closing price $y_0$ (Yuan) | The last trading price of Shanghai Composite Index on a certain day |
| Amplitude $y_1$ % | The percentage of the absolute value of difference between the highest price and the lowest price on a given day after opening, relative to the closing price that day; it reflects the activity degree of stocks to some extent. |
| Total lots $y_2$ (10000) | The total trade volume that day. |
| Amount $y_3$ (0.1 billion) | Total trading amount determined according to the total shares of a stock and the corresponding trading prices. |
| DEA $y_4$ | smoothed difference exponential average curve, the difference between one long-term and one short-term smoothed average curves, and used to judge market trading |
| Price movement $y_5$ (Yuan) | The difference between the closing price on a given day and the closing price on the previous day, and used to judge whether the stock price rises or declines. |
| DIFF $y_6$ | The difference between the short-term and long-term index smoothing moving average lines of the closing price, and often combined with DEA to judge the stock trading situation |
| Rising amplitude $y_7$ % | The percentage of rising in the closing price that day relative to the closing price on the previous day. |
| 5-day average $y_8$ | The average value of trading amounts on 5 days, and reflecting the trend of stock trading amounts |

To eliminate differences in dimensions, the data were standardized. The time sequences of the closing price of Shanghai Composite Index after the standardization were shown below.

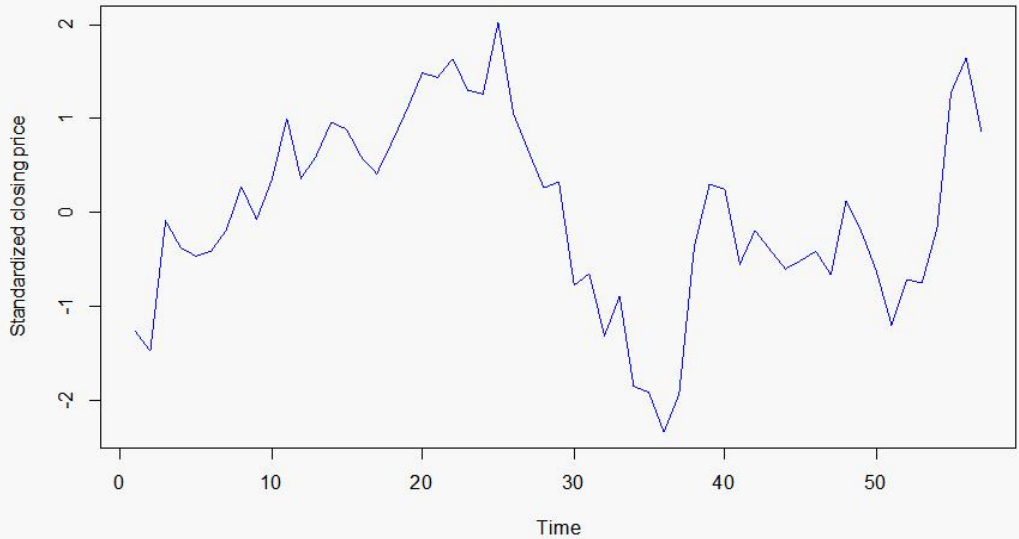

**Fig. 4** The closing time diagram of the closing price

Fig. 4 shows the trend of the closing price of Shanghai Composite Index. Clearly, after a slump at early 2016, the closing price generally was rising until March 2018, and fluctuated to different extents in this period. Since March 2018 when trade disputes between China and the United States were aggravated, the closing price of Shanghai Composite Index was slumping until middle February 2019. In recent two years, the trend was still unclear. At early March 2020, the Shanghai Composite Index minimized again. On the one hand, the unlimited quantitative easing from the US Federal Reserve indeed stabilized the stock market. The financial relief policies from the U.S. also helped relieve the economic impact due to the pandemic and oil price slip. On the other hand, the Interim Measures for the Administration of Insurance Asset Management Products released by China Banking and Insurance Regulatory Commission contributed to guide long-term participation of funds into the capital market. The minister of the Chinese Ministry of Finance stated that China's active financial policies will be further enthusiastic and promising, and the excitement policies from the Chinese Government will be strengthened. Thus, the Shanghai Composite Index gradually rises and the global stock indices are rebounding persistently. For these reasons, the monthly data of Shanghai Composite Index were used to build a differential equation model of closing price.

**4.2 Analysis of Shanghai Composite Index based on dynamic differential equation**

Let the step length be $h=2$ in the monthly data of Shanghai Index 300 from January 2016 to the end of October 2020. Then the first-order differences in adjacent closing prices of Shanghai Composite Index were computed using Eq. (3). Next, the data of independent variables were iteratively computed as stated above. The trend of first-order differences from the standardized closing prices $\Delta y_t$ was shown in Fig. 5.

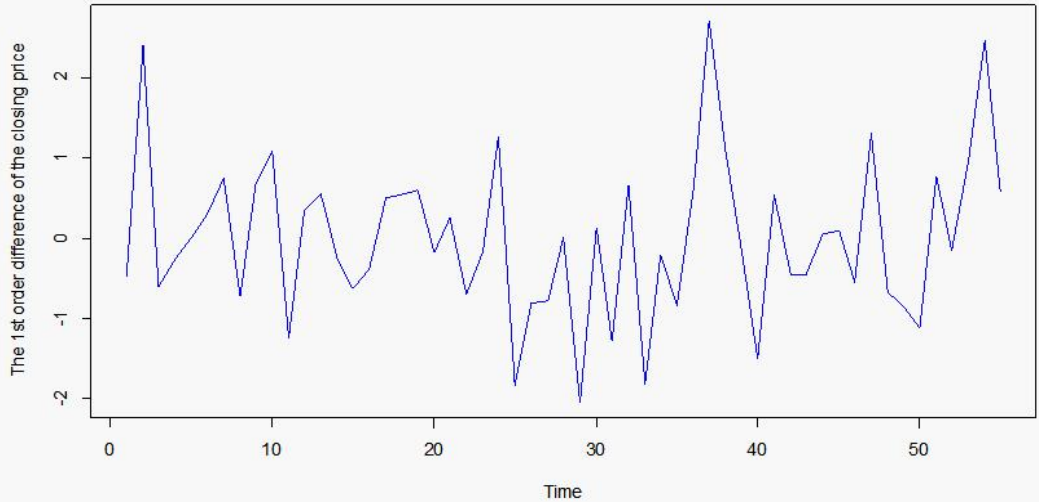

**Fig. 5** First-order difference of closing price

As shown in Fig. 5, the first-order differences of closing prices change very frequently, and can be modeled and analyzed using a linear regression model. Then the data were used to build multiple linear regression models. Model averaging was conducted using the Optimal weight selection method.

In a differential equation, according to the fitting results of the models and the significance analysis of variables, we chose $y_1, y_4, y_6$ as the main regression variable and other variables as auxiliary variables. Seven linear models were chosen as the submodels of the Optimal model averaging process.

**Table 8** Model output selection

|  | $y_1$ | $y_2$ | $y_3$ | $y_4$ | $y_5$ | $y_6$ | $y_7$ | $y_8$ | weight $\lambda_i$ |
|---|---|---|---|---|---|---|---|---|---|
| Model 1 | ✓ | ✓ | ✓ | ✓ | ✓ | ✓ | ✓ | ✓ | 0.8541 |
| Model 2 | ✓ | ✓ | ✓ | ✓ | ✓ | ✓ | ✓ |  | 0.1046 |
| Model 3 | ✓ | ✓ | ✓ | ✓ | ✓ | ✓ |  |  | 0.0132 |
| Model 4 | ✓ |  |  | ✓ |  | ✓ |  |  | 2.27e-05 |
| Model 5 | ✓ | ✓ | ✓ | ✓ |  | ✓ |  |  | 0.0017 |
| Model 6 | ✓ | ✓ | ✓ | ✓ |  | ✓ | ✓ |  | 0.0132 |
| Model 7 | ✓ | ✓ | ✓ | ✓ |  | ✓ |  | ✓ | 0.0134 |

The weighted averages of the above models were determined. Then the differential equation averaging model was obtained.

This model indicates the first-order differences of closing prices are positively correlated with the amplitude, total amount, price movement, rising amplitude, and 5-day average during the R-K process, and are negatively correlated with total lots, DEA and DIFF. Then the MSEs of the model averaging were compared with the MSEs of other submodels.

**Table 9** Comparison of models

| Models | Standard error | Mean square error | Mean absolute error | Goodness of fit |
|---|---|---|---|---|
| Model 1 | 0.5836 | 0.3406 | 0.0786 | 0.6539 |
| Model 2 | 0.5888 | 0.3467 | 0.0793 | 0.6468 |
| Model 3 | 0.5893 | 0.3473 | 0.0794 | 0.6463 |
| Model 4 | 0.6109 | 0.3732 | 0.0823 | 0.6199 |
| Model 5 | 0.5903 | 0.3485 | 0.0796 | 0.6450 |
| Model 6 | 0.5892 | 0.3471 | 0.0794 | 0.6464 |
| Model 7 | 0.5854 | 0.3427 | 0.0789 | 0.6510 |
| FMA model | 0.5836 | 0.3407 | 0.0787 | 0.6536 |

The MSEs of the differential averaging model are smaller than those of most submodels (Table 9). The aim of Optimal weight selection is to minimize the MSEs of models. Thus, the fitting results of the differential averaging model are more stable than most other models. The fitting results are shown in Fig. 6.

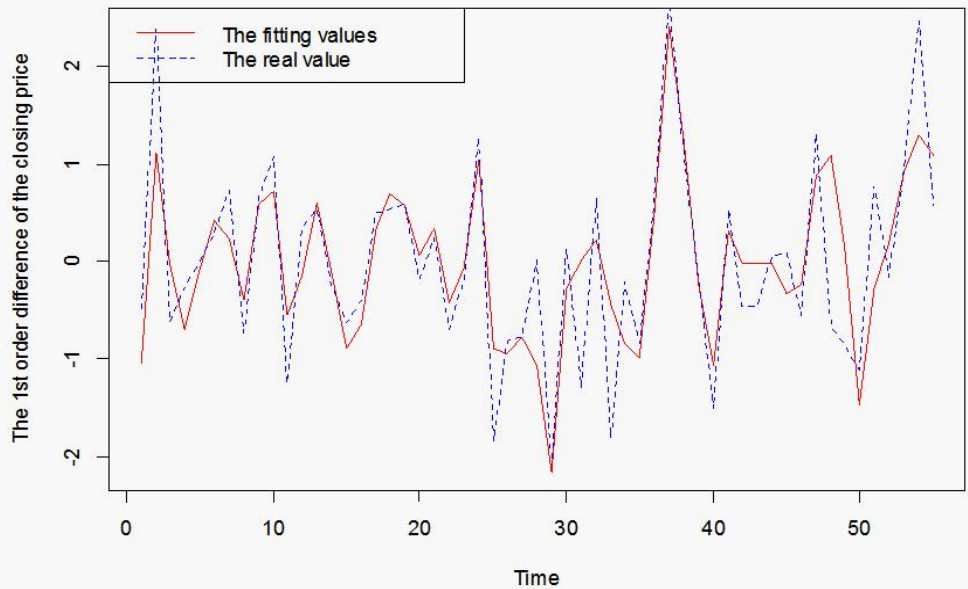

**Fig. 6** Rendering of first-order differential fitting of closing price

The linear averaging models can well fit the part of closing prices with large variation in first-order difference, but the fitting trend of closing prices with small variation in first-order differences is locally different from the real values (Fig. 6). Generally, the first-order differences from the model fitting are close to the real differences. Then the residual errors of this model were sent to normality test. The histograms and QQ maps of residual errors were shown in Fig. 7.

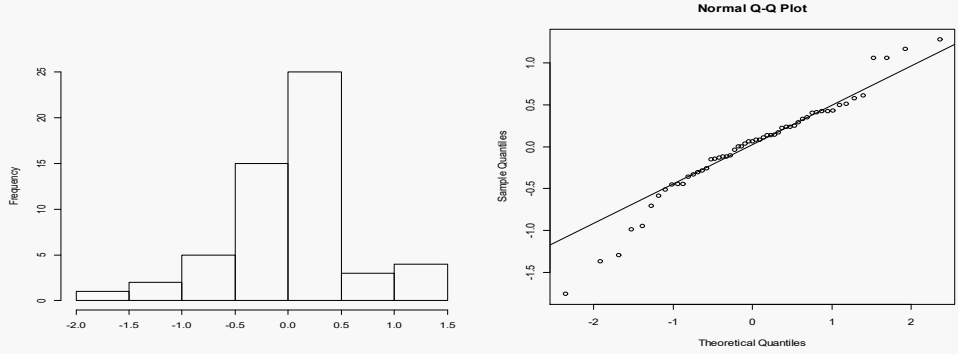

**Fig. 7** Distribution of model residual errors

The histograms in Fig. 7 show that the residual errors of this model are concentrated at 0. The majority of sample points on the QQ maps are approximately distributed near the straight line, indicating the residual errors basically obey normal distribution and accord with the hypotheses.

Next, the fitting values $\hat{y}_{t+1}$ were determined by adding the data $\Delta y_t$ at each time point from the regression model with the closing price $y_t$ at the corresponding time. After that, an ordinary linear averaging model was built to compare the fitting results.

The above data were also used in the linear regression model. The regression variables after Runge-Kutta dispersion largely decreased intervariable correlations, but the intervariable independence shall be maintained maximally during linear modeling. Thus, a least square model, a ridge regression model and a generalized linear model were built, and the heteroscedasticity was reduced using Optimal weight selection to make the models more stable.

Analysis with the linear models shows that a larger amplitude of $y_1$ means the stock activeness is higher, and the closing price on that day declines. Moreover, the total trading volume $y_2$ is larger and the closing price is lower. Generally, the Chinese stock market involves numerous retail investors. When the stock prices drop, shareholders take this opportunity to buy stocks and accordingly, the stock activeness is increased. The larger trading amount $y_3$ means the closing price is larger, so more shareholders sell out stocks. $y_4 \sim y_8$ belong to auxiliary variables. Similarly, the estimation coefficient of trading amount is also consistent with the dynamic differential equation model, which will not be elaborated here.

Next, the fitting results from the differential averaging model and the linear averaging models were compared (Table 10).

**Table 10** Comparison of differential averaging model and linear averaging models

|  | Standard error | Mean square error | Mean absolute error | Goodness of fit |
|---|---|---|---|---|
| Average model of differential equations | 0.5836 | 0.3407 | 0.0787 | 0.6536 |
| Linear average model | 0.7383 | 0.4879 | 0.0925 | 0.4900 |

The indices in the differential averaging model are smaller than those of the linear averaging models (Table 10) and especially, the MSEs are smaller than those of the static linear model. Then the MSEs of the two models were compared.

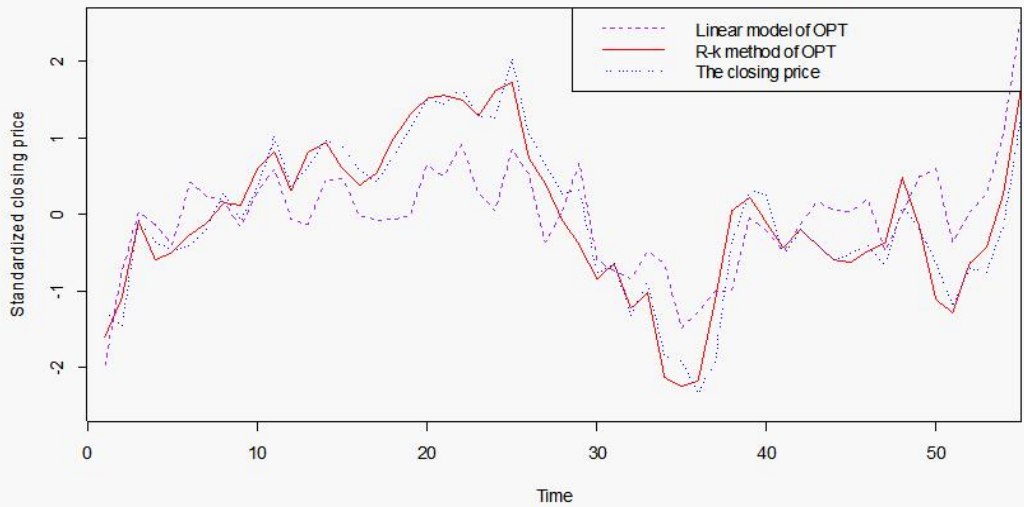

**Fig. 8** Comparison of model fitting

Clearly, the regression model derived from the dynamic differential equation better fits the closing prices compared with ordinary linear models, and basically reflects the value and trend of closing price at each period. The ordinary linear models can only roughly reflect the trend of closing price, and only a small part of its fitted values are close to the real values. Moreover, after the regression variables were discretized via the Runge-Kutta method, the significance of intervariable correlations decreased. Thus, the differential averaging model is more precise than ordinary linear averaging models. This method integrates the fluctuation features of multiple Shanghai Composite Index fitting models, and provides a methodological reference for analyzing stock market situations. This method also offers a new clue for regression analysis in statistics. In other words, when the analytical results of ordinary linear regression models are non-ideal, a differential model can be built.

## 5. Conclusions

A differential model was built, and then the continuous differential model was discretized using the classic four-order Runge-Kutta method. Then the parameters of the discretized model were estimated by combining the model averaging Optimal weight selection method. The consistency of the estimated model parameters was validated. Then data of Shanghai Composite Index were analyzed, and regression models were built. Next, these results were compared with the ordinary linear regression models with Optimal weight selection. The differential averaging model with RK4 can yield better results than ordinary linear averaging models. In practical applications, before model discretization, we shall check whether the variables are continuous prior to the setting of step length. Additionally, this method is only applicable to problems with time series. In the future, we will expand this method to more models. This method was used to analyzed Shanghai Composite Index from the perspective of differential equations. Compared with traditional linear regression models, this method provides new clues and schemes for research with differential models and other linear models. This method can also be expanded to other meaningful fields, such as parameter estimation of nonlinear autonomous differential equations, which is another direction of our future research.

**Acknowledgments** This work was supported by the National Science Foundation of Jilin Province under Grant No.20200201273JC.